\documentclass[11pt]{amsart}
\usepackage{amsmath}
\usepackage{mathtools}
\usepackage{amssymb}
\usepackage{fullpage}
\usepackage{mathrsfs}
\usepackage{mathptmx}
\usepackage{enumerate}
\usepackage{stackrel}
\usepackage{xcolor}
\usepackage{hyperref}

\usepackage{amsopn}
\usepackage{dsfont}
\newcommand{\hide}[1]{}
\usepackage{enumerate}
\numberwithin{equation}{section}
\usepackage{float}
\usepackage{graphicx}
\usepackage{caption}
\usepackage{subcaption}
\usepackage{titlesec}
\usepackage{fancyhdr}
\usepackage{textcomp}
\usepackage{fancyhdr}
\newtheorem{corollary}{Corollary}[section]
\newtheorem{lemma}{Lemma}[section]
\newtheorem{theorem}{Theorem}[section]
\newtheorem{proposition}{Proposition}[section]
\newtheorem{definition}{Definition}[section]

\theoremstyle{definition}
\newtheorem{remark}{Remark}[section]

\DeclareMathOperator{\Log}{Log}

\DeclareMathOperator{\dist}{dist}
\DeclareMathOperator{\arctanh}{arctanh}

\DeclareMathOperator{\Arg}{Arg}

\DeclareMathOperator{\Spir}{Spir}

\DeclareMathOperator{\di}{d}
\DeclareMathOperator{\D}{\mathbb{D}}
\DeclareMathOperator{\RE}{Re}
\DeclareMathOperator{\IM}{Im}
\titleformat{\section}
  {\normalfont\scshape}{\thesection}{1em}{\centering}
\titleformat{\subsection}[runin]
  {\bfseries}{\thesubsection}{1em}{}

\begin{document}

\title{Speeds of Convergence for Petals of Semigroups of Holomorphic Functions}

\author{Maria Kourou$^{\S}$}
\thanks{$^{\S}$Partially supported by the Alexander von Humboldt Foundation.}  
\address{Department of Mathematics, Julius-Maximilians University of Wuerzburg, 97074, Wuerzburg, Germany}
\email{maria.kourou@mathematik.uni-wuerzburg.de}   

\author{Konstantinos Zarvalis}  
\address{Department of Mathematics, Aristotle University of Thessaloniki, 54124, Thessaloniki, Greece}
\email{zarkonath@math.auth.gr}   

\fancyhf{}
\renewcommand{\headrulewidth}{0.4pt}
\fancyhead[RO,LE]{\small \thepage}
\fancyhead[CE]{\footnotesize M. KOUROU \& K. ZARVALIS}
\fancyhead[CO]{\footnotesize SPEEDS OF CONVERGENCE FOR PETALS OF SEMIGROUPS OF HOLOMORPHIC FUNCTIONS}
\fancyfoot[L,R,C]{}

\date{}

\begin{abstract}
We study the backward dynamics of one-parameter semigroups of holomorphic self-maps of the unit disk. More specifically, we introduce the speeds of convergence for petals of the semigroup, namely the total, orthogonal, and tangential speeds. These are analogous to speeds of convergence introduced by Bracci, yet profoundly different due to the nature of backward dynamics.  
Results are extracted on the asymptotic behavior of speeds of petals, depending on the type of the petal. 
We further discuss the asymptotic behavior of the hyperbolic distance along non-regular backward orbits. 
\end{abstract}

\maketitle

\section{Introduction}

One-parameter continuous semigroups of holomorphic functions in the unit disk, or from now on \textit{semigroups in} $\mathbb{D}$, have stimulated the scientific interest in recent years. 
The introduction to their present form was made by Berkson and Porta in \cite{berksonporta}, as a direct aspect of semigroups of composition operators. Later, Contreras and D{\'i}az-Madrigal \cite{AnalyticFlows} established their main characteristics leading to a plethora of new results. 
A thorough analysis as well as recent advances on semigroups in $\mathbb{D}$ can be found in the recent monograph \cite{Booksem} and references therein.

A semigroup in $\mathbb{D}$ is a family $(\phi_t)_{t\ge0}$ of holomorphic self-maps of the unit disk that satisfy the following conditions:
\begin{enumerate}[(i)]
\item $\phi_0(z)=z$, for all $z \in \D$;
\item $\phi_{t+s}(z) = \phi_t \left( \phi_s (z) \right)$, for every $t,s \geq 0$ and $z \in \D$;
\item $\phi_t(z) \xrightarrow{t \to 0^{+}} z$, uniformly on compacta in $\D$.
\end{enumerate}

 If, in addition, for some (equivalently all) $t_0 > 0$ it is true that $\phi_{t_0}$ is an automorphism of $\D$, then $(\phi_t)$ is called a \textit{group}. 
For all semigroups that are not groups, the continuous version of the Denjoy--Wolff Theorem asserts the existence of a unique point $\tau \in \overline{\D}$ such that $\phi_t(z) \to\tau$, as $t\to +\infty$, for all $z\in\mathbb{D}.$ 
This point $\tau$ is called the \textit{Denjoy--Wolff point} of the semigroup; see \cite[Theorem 1.4.17]{abate}. If $\tau\in\mathbb{D}$, the semigroup is characterized as \textit{elliptic}, while if the Denjoy--Wolff point lies on the unit circle, $(\phi_t)$ is called \textit{non-elliptic}.

Both for elliptic and non-elliptic semigroups $(\phi_t)$, there exists some complex number $\mu$ with $ \RE \mu>0$ such that the (angular) derivative  $\phi_t^{'}(\tau)=e^{-\mu t}$, for all $t\ge0$. Note that in the case of non-elliptic semigroups, $\mu\geq 0$.
The number $\mu$ is called the \textit{spectral value} of the semigroup. 
 
Fix $z \in \D$. The curve $ \gamma_z: [0, + \infty)  \to \D$ with $\gamma_z(t)= \phi_t(z)$ is called the \textit{trajectory} of $z$. Clearly, for all points $z\in \D$ the trajectory converges to $\tau$, as $t\to +\infty$. 
Bracci \cite{brspeeds} introduced three novel quantities concerning trajectories of non-elliptic semigroups of holomorphic functions. These quantities, the so-called \textit{speeds}, provide interesting results with regard to the rate of convergence of trajectories to the Denjoy--Wolff point of the semigroup. Let $(\phi_t)$ be a non-elliptic semigroup in the unit disk $\mathbb{D}$ with Denjoy--Wolff point $\tau$. Consider $\gamma:(-1,1)\to\mathbb{D}$, with $\gamma(r)=r\tau$, to be the diameter of the unit disk with endpoints $\tau$ and $-\tau$. Clearly, $\gamma$ is a geodesic for the hyperbolic distance $d_\mathbb{D}$ of the unit disk. We denote by $\pi_\gamma(z)$ the \textit{projection} of a point $z\in\mathbb{D}$ onto the curve $\gamma$, which satisfies
$$d_\mathbb{D}(z,\pi_\gamma(z))=d_\mathbb{D}(z,\gamma):=\inf\limits_{r\in(-1,1)}d_\mathbb{D}(z,\gamma(r)).$$
The function 
$$v(t)=d_\mathbb{D}(0,\phi_t(0)), \quad t\ge0,$$
is called \textit{total speed} of $(\phi_t)$. This can be decomposed into two other functions, the \textit{orthogonal speed} of $(\phi_t)$
$$v^o(t)=d_\mathbb{D}(0,\pi_\gamma(\phi_t(0))), \quad t\ge0,$$ and the \textit{tangential speed} of $(\phi_t)$
$$v^T(t)=d_\mathbb{D}(\phi_t(0),\pi_\gamma(\phi_t(0)))=d_\mathbb{D}(\phi_t(0),\gamma), \quad t\ge0.$$

With a first glance, it seems as if the three speeds are solely defined with respect to the trajectory with starting point $0$. Nevertheless, Bracci proved that asymptotically these functions do not depend on the starting point and therefore the selection of any point $z\in\mathbb{D}$ instead of $0$ is eligible.
 Furthermore, observing the asymptotic behavior of the speeds, we obtain the type of the semigroup. According to \cite[Proposition 6.1]{brspeeds}, 
\begin{equation}\label{speeds2}
 \lim_{t \to +\infty} \frac{v(t)}{t}= \lim_{t \to +\infty} \frac{v^{o}(t)}{t}= \frac{\mu}{2},
\end{equation}
where $\mu$ denotes the spectral value of $(\phi_t)$, a result that agrees with that in \cite{arosiobracci} and \cite[Lemma 9.1.2, Theorem 9.1.9]{Booksem}.

The definition of speeds of convergence and the statement of a variety of questions in \cite{brspeeds} ignited the research interest and was the stepping stone for several works. 
The second named author disproved in \cite{Zarvalis} a conjecture on the upper bound for the tangential speed in parabolic semigroups. Cordella in \cite{Cordella} worked on the asymptotic upper bound for the tangential speed. 
Concerning the orthogonal speed, Bracci, Cordella, and the first named author \cite{brcorkou} examined its asymptotic monotonicity with respect to semigroups for a variety of cases. 
Quite recently, Betsakos and Karamanlis \cite{BetKar2022} generalized the above result and proved by means of harmonic measure the monotonicity of the orthogonal speed for all cases of non-elliptic semigroups. 

The main focus of the current article is to establish speeds of convergence in the setting of the backward dynamics of a semigroup of holomorphic self-maps of $\D$. 
Advancements in the direction of backward dynamics are made by Elin, Shoikhet, and Zalcman \cite{ESZ2008} as well as by Bracci, Contreras, D\'iaz-Madrigal, and Gaussier \cite{petals}, where the theory on the backward flow of a semigroup is concretely settled. 
The \textit{backward invariant set} of a semigroup $(\phi_t)$ of holomorphic self-maps of $\D$ is defined as $$\mathcal{W} : =\bigcap_{t\geq 0} \phi_t(\D)$$ and it is exactly the set where the restriction of the semigroup is a group of automorphisms. Every non-empty connected component of the interior of $\mathcal{W}$ is called a \textit{petal} of $(\phi_t)$. 
In the course of the paper, we work with semigroups, which are not groups and whose backward invariant set is non-empty. Thus, we exclude the trivial cases from a backward-dynamical point of view of $\mathcal{W}$ being empty or equal to the whole unit disk $\D$. 

As already mentioned, \emph{forward} speeds is an attribute of the semigroup and their asymptotic behavior is similar regardless of the chosen initial point $z\in \D$, due to the attractive nature of the Denjoy--Wolff point. However, as emphasized in the theory of backward dynamics (see \cite[Chapter 13]{Booksem} or Section \ref{semigroups}), this is not the case in the backward-dynamical setting. 
Due to the lack of an analogous Denjoy-Wolff Theorem, speeds of convergence can only be defined for points in the backward invariant set of a semigroup. 
As it is clear in Propositions \ref{prop:total} and \ref{prop:ort&tang}, the need to restrict to the geometry of a petal of a semigroup arises in order to examine the behavior of the so-called \textit{backward speeds}. Definition \ref{def:speeds} establishes the expression for the \textit{total speed of a petal} and further illustrates the suitable geodesic selection that leads to the definition of the \textit{orthogonal} and \textit{tangential speeds of a petal}.  

The asymptotic behavior of speeds of petals are outlined in Sections \ref{sec:totalspeed} and \ref{sec:orthogonal&tandspeed}. We prove that the asymptotic behavior of the total speed of a petal solely depends on the type of the petal and its geometry, as expressed in detail in Theorems \ref{thm:divratehyperbolic} and \ref{thm:bounds}. 
Concerning the tangential speed, we examine under which circumstances its limit superior is finite. As stated in Theorem \ref{thm:asympttang}, the finiteness of the limit depends on the type of the petal. 
A connection between all speeds of a certain petal is provided by a generalization of Bracci's Pythagoras Theorem. Taking all these into account, we obtain information on the asymptotic behavior of the orthogonal speed of the petal; see Corollary \ref{cor:orthogonalasymb}. 

Section \ref{sec:non-regular} is devoted to a special case of backward dynamics of a semigroup of holomorphic self-maps of $\D$. We extend the definition of the total speed for points that lie on the boundary of the backward invariant set, thus they do not lie inside any petal. 
The initiative for this study is the better understanding of the nature of the hyperbolic distance of $\D$ along the so-called \textit{non-regular backward orbits}; see Section \ref{semigroups} for precise definitions. As observed in Proposition \ref{prop:non-regular}, the asymptotic behavior of the total speed in this case depends on the geometry of the associated Koenigs domain of the semigroup and on the intrinsic characteristics of the semigroup. 

The proofs presented in the current work rely on conformally invariant tools such as harmonic measure, extremal length and hyperbolic geometry. Meanwhile, some basic characteristics of semigroups are further utilized such as the geometry of the associated Koenigs domain and the associated infinitesimal generator. A collection of the background information on the above subjects follows in Section \ref{background}.

\section{Preparation for the Proofs}\label{background}

\subsection{ Hyperbolic Metric }\label{hypmetric}

We briefly present a variety of conformally invariant quantities that will be of use throughout the course of the proofs. 
We start with some information about hyperbolic geometry (see \cite[Chapter 5]{Booksem}). 
The \textit{density} of the \textit{hyperbolic metric} in $\D$ is $ \lambda_{\D} (z)= (1-|z|^2)^{-1}$.

The \textit{hyperbolic distance} between two points $z,w$ in the unit disk is 
\begin{equation}\label{hyperbolicdistance}
d_{\D} (z,w)= \arctanh \rho_{\D}(z,w), \quad \text{where} \quad \rho_{\D}(z,w) = \left| \frac{z-w}{1- \bar{z}w} \right| 
\end{equation}
 denotes the \textit{pseudo-hyperbolic distance} in the unit disk $\D$. 

Using conformal mappings, we are able to transcend the notions of hyperbolic metric and distance to any simply connected domain, other than the complex plane. Indeed, let $\Omega\subsetneq\mathbb{C}$ be a simply connected domain and let $f:\Omega\to\mathbb{D}$ be a Riemann map. Then the hyperbolic distance between two points $z,w$ in $\Omega$ is given by
$d_\Omega(z,w):=d_\mathbb{D}(f(z),f(w))$
and hence the hyperbolic distance, which does not depend on the choice of the Riemann map $f$, is a conformally invariant quantity. 

Applying known conformal mappings, we can calculate the hyperbolic distance in certain domains. For example, if $S:=\{ z \in \mathbb{C} : |\IM z| < \frac{\pi}{2} \}$ with $z,w\in S$, and $\mathbb{H}:=\{ w \in \mathbb{C} : \RE w >0 \}$ with $x,y \in \mathbb{H}$, then the hyperbolic distances in the respective domains are given by the formula
\begin{equation}\label{eq:hypdiststripandhalfplane} 
    d_S(z,w)  = \arctanh \left| \frac{e^z - e^w }{e^z + e^{\overline{w}}} \right| \quad \text{and} \quad d_{\mathbb{H}} (x,y) = \arctanh \left| \frac{x-y}{x+\overline{y}} \right|.
\end{equation}

Other than the conformal invariance, a very important property of the hyperbolic distance is its domain monotonicity. To be exact, let $\Omega_1,\Omega_2\subsetneq\mathbb{C}$ be two simply connected domains with $\Omega_1\subset\Omega_2$. Then $d_{\Omega_1}(z,w)\ge d_{\Omega_2}(z,w)$, for all $z,w\in\Omega_1$. 
 
Let $\Omega\subsetneq\mathbb{C}$ be a simply connected domain and $\delta_\Omega(z):=\dist(z,\partial \Omega)$ denote the Euclidean distance of $z$ from the boundary of $\Omega$. The hyperbolic density $\lambda_\Omega$ shares a deep connection with $\delta_{\Omega}$, as observed in the following \textit{Distance Lemma}.
\begin{lemma}{\cite[Theorem 5.3.1]{Booksem}}\label{lem:distancelemma}
    Let $\Omega\subsetneq\mathbb{C}$ be a simply connected domain. Then, for all $z,w\in \Omega$,
    $$\frac{1}{4}\log\left(1+\frac{|z-w|}{\min\{\delta_\Omega(z),\delta_\Omega(w)\}}\right)\le d_\Omega(z,w)\le\int\limits_\Gamma\frac{|d\zeta|}{\delta_\Omega(\zeta)},$$
    where $\Gamma$ can be any piecewise $C^1$-smooth curve joining $z$ and $w$ inside $\Omega$.
\end{lemma}

\subsection{ Harmonic Measure}\label{harmonicmeasure}

 Let $\Omega$ be a proper subdomain of $\widehat{\mathbb{C}}$ with non-polar boundary and $\mathcal{B}(\partial\Omega)$ be the $\sigma$-algebra of all Borel sets of $\partial\Omega$.   
Suppose $E \in \mathcal{B}(\partial\Omega)$. The harmonic measure of $E$ with respect to $\Omega$ at a point $z \in\Omega$ is the solution of the generalized Dirichlet problem for the Laplacian in $\Omega$ with boundary values $1$ on $E$ and $0$ on $\partial\Omega\setminus E$. 

For a fixed $E \in \mathcal{B}(\partial \Omega)$, $\omega ( \cdot , E, \Omega)$ is a harmonic and bounded function on $\Omega$. 
In addition, for a fixed point $z \in \Omega$, the map 
$ \omega ( z, \cdot, \Omega) : \mathcal{B}(\partial \Omega)  \to [0,1]$, with $E \mapsto \omega ( z, E , \Omega)$,
is a Borel probability measure on $\partial \Omega$. 

A major property of the harmonic measure is its conformal invariance; see \cite[\S\ 4.3]{ransford}. In addition, the harmonic measure conceals a useful monotonicity property. Let $\Omega_1,\Omega_2\subsetneq\widehat{\mathbb{C}}$ be two simply connected domains with non-polar boundary satisfying $\Omega_1\subset\Omega_2$ and let $E\subset\partial\Omega_1\cap\partial\Omega_2$ be a Borel set. Then
$\omega(z,E,\Omega_1)\le\omega(z,E,\Omega_2)$, for all $z\in\Omega_1$.

The aforementioned monotonicity property can be made more precise by means of the following Proposition.
\begin{proposition}[Strong Markov Property for Harmonic Measure]\cite[p.88]{portstone}\label{markov}
Let $\Omega_1,\Omega_2$ and $E$ be as above. Then, for all $z\in\Omega_1$, we have 
$$\omega(z,E,\Omega_2) =\omega(z,E,\Omega_1)+\int\limits_{\partial\Omega_1\setminus\partial\Omega_2} \omega (\zeta, E ,\Omega_2) \cdot \omega(z,d\zeta, \Omega_1). $$
\end{proposition} 

\begin{remark}\label{caratheodory}
    Let $E\subset\partial\mathbb{D}$ be a circular arc with endpoints $a$ and $b$. Then, the level set
    $$L_k=\{z\in\mathbb{D}:\omega(z,E,\mathbb{D})=k\}, \quad 0<k<1,$$
    is a circular arc (or a diameter if $E$ is a half-circle and $k=\frac{1}{2}$) with endpoints $a$ and $b$ that meets the unit circle with angle $k\pi$ (see \cite[p. 155]{carathe}). As a result, a sequence $\{z_n\}\subset\mathbb{D}$ satisfying $\lim\limits_{n\to+\infty}z_n=a$ (or $b$) converges by angle $\theta\in[0,\pi]$ if and only if $\lim\limits_{n\to+\infty}\omega(z_n,E,\mathbb{D})=\frac{\theta}{\pi}$. In particular, the sequence converges tangentially if and only if the above limit is equal to $0$ or $1$. Naturally, this result also holds for continuous curves.
\end{remark}

More information and detailed theory of harmonic measure can be found in \cite{ahlforsconformal, margarnett} and \cite[Chapter 4]{ransford}. 

\subsection{Extremal Length}

Let $\Omega\subsetneq\mathbb{C}$ be a domain and suppose that $\Gamma$ is a family of locally rectifiable curves inside $\Omega$. For a non-negative Borel measurable in $\Omega$ function $\rho$, we define the $\rho$\textit{-area} of $\Omega$
$$A(\rho,\Omega):=\iint\limits_\Omega\rho(z)^2dxdy, \quad z=x+iy.$$
We also define the $\rho$\textit{-length} of $\Gamma$ as
$$L(\rho,\Gamma):=\inf\limits_{\gamma\in\Gamma}\int\limits_\gamma\rho(z)|dz|.$$
Then, the extremal length with respect to $\Omega$ of $\Gamma$ is the quantity
$$\lambda_\Omega(\Gamma)=\sup\limits_\rho\frac{L(\rho,\Gamma)^2}{A(\rho,\Omega)},$$
where the supremum is taken over all non-negative Borel measurable functions $\rho$ whose $\rho$-area is a positive real number.

Let us now suppose that $\Gamma_1,\Gamma_2$ are two distinct families of locally rectifiable curves defined on two domains $\Omega_1,\Omega_2\subsetneq\mathbb{C}$, respectively. The very definition of extremal length provides a type of a monotonicity property: if $\Gamma_1$ is a subfamily of $\Gamma_2$, then $\lambda_{\Omega_1}(\Gamma_1)\ge\lambda_{\Omega_2}(\Gamma_2)$. A second type of a monotonicity property derives from the following:
\begin{lemma}{\cite[The Extension Rule, \S IV.3]{margarnett}}\label{extr-length-monoton}
    If every curve belonging to the family $\Gamma_1$ contains a subcurve that belongs to the family $\Gamma_2$, then $\lambda_{\Omega_1}(\Gamma_1)\ge\lambda_{\Omega_2}(\Gamma_2)$. 
\end{lemma}
Along with the two aforementioned monotonicity properties, we use one significant instance of the extremal length. Once again, let $\Omega\subsetneq\mathbb{C}$ be a domain and $E,F\subset\overline{\Omega}$. The \textit{extremal distance} between $E$ and $F$ inside $\Omega$ is defined as the quantity
$$\lambda(E,F,\Omega):=\lambda_{\Omega}(\Gamma),$$
where $\Gamma$ is the family of all locally rectifiable curves joining $E$ to $F$ inside $\Omega$. Trivially, the extremal distance is also conformally invariant by extension, due to the fact that the extremal length is conformally invariant. 

A useful tool that connects extremal length and harmonic measure is the so-called \textit{Beurling's Estimate}.
\begin{lemma}{\cite[Theorem IV.5.3]{margarnett}}\label{beurling}
Let $\Omega$ be a a simply connected Jordan domain and let $E$ be a finite union of arcs lying on $\partial\Omega$. Let $z\in\Omega$. Then,
$$\omega(z,E,\Omega)\le C\cdot e^{-\pi\cdot\lambda(\gamma,E,\Omega)}$$
for every Jordan arc $\gamma$ joining $z$ to $\partial\Omega\setminus E$ inside $\Omega$, where $C>0$ is an absolute constant.
\end{lemma}
For a thorough presentation of the theory of extremal length, the interested reader may refer to \cite[Chapter 4]{ahlforsconformal} and \cite[Chapter IV]{margarnett}.

\subsection{ One-Parameter Semigroups - Koenigs Function}
As stated in the Introduction, we exclude one-parameter groups from the spectrum of our study, due to their trivial behavior in the backward-dynamical setting. Hence, from this point on, we refer to one-parameter semigroups that are not groups, as simply \textit{one-parameter semigroups}.

Non-elliptic semigroups are further divided into two categories based on their spectral value $\mu$. For $\mu>0$, a semigroup is characterized as \textit{hyperbolic}, whereas, for $\mu=0$, the semigroup is called \textit{parabolic}.
For every one-parameter non-elliptic semigroup, there exists a unique Riemann mapping $h$ fixing the origin such that 
$$h(\phi_t(z))=h(z)+t, \quad \text{for all }z\in\mathbb{D} \text{ and }t\ge0.$$
The simply connected domain $\Omega:=h(\mathbb{D})$ is called the \textit{Koenigs domain} (also known as \textit{associated planar domain}) of the semigroup. The Koenigs domain $\Omega$ of a non-elliptic semigroup is \textit{convex in the positive direction}, which means that $\{ w+t: t\ge0 \} \subset \Omega$, for every $w \in \Omega$. By its definition, it can be understood that a major benefit of the Koenigs function is the linearization of the trajectories of the points in $\D$ under the semigroup.  Basically, the image through $h$ of the trajectory of some $z\in\mathbb{D}$ is a horizontal half-line emanating from $h(z)$ towards $\infty$ in the positive direction (i.e. with constant imaginary part and increasing real part).

The Koenigs function can also be uniquely defined for one-parameter elliptic semigroups. In this case, $h(\tau) =0$, where $\tau \in \D$ is the Denjoy--Wolff point of the semigroup $(\phi_t)$ and 
\begin{equation}\label{mapofhelliptic}
h(\phi_t(z)) = e^{-\mu t} h(z),\quad \text{for all }  z\in \D \text{ and } t\ge0, 
\end{equation}
 where $\mu$ is the spectral value of $(\phi_t)$. This time, the Koenigs domain $\Omega$ is \textit{$\mu$-spirallike} with respect to $0$, since $0\in\Omega$ and $e^{-\mu t }\Omega \subseteq \Omega$, for all $t \geq 0$. 
\hide{In addition, every point $w\in\mathbb{C}\setminus \{0\}$ can be written in the form $w=e^{-\mu t+ i \theta}$, for some $t \in \mathbb{R}$ and $\theta\in [-\pi, \pi)$. 
The \textit{$\mu$-spirallike argument} of $w$ is defined as $\Arg_{\mu}(w):=\theta.$}
The Koenigs function associated to an elliptic semigroup maps the trajectory of a $z\in\mathbb{D}$ onto the half-spiral, 
which winds around $h(\tau)=0$ infinitely many times, as $t\to+\infty$. 


Further information on the Koenigs function of a one-parameter semigroup and advances on its geometry and overall behavior can be found in \cite[Chapter 9]{Booksem} as well as in references therein.


\subsection{Backward Orbits}\label{semigroups}

As further stated in the Introduction, every trajectory of a semigroup converges to the Denjoy--Wolff point of the semigroup. However, the asymptotic behavior of backward orbits is not so straightforward. Before getting into detail about the convergence of backward orbits, we need some important notions.
\begin{definition}\cite[Chapter 12]{Booksem}
\normalfont
Let $(\phi_t)$ be a semigroup in $\D$. A point $\sigma\in\partial\mathbb{D}$ is called a \textit{boundary fixed point} of $(\phi_t)$ provided $\angle\lim\limits_{z\to\sigma}\phi_t(z)=\sigma$, for all $t\ge0$. A boundary fixed point of $(\phi_t)$ is characterized as \emph{regular}, if the angular derivative of $\phi_t$ at $\sigma$ is finite, for all $t \geq 0$. Otherwise, we say that it is \textit{non-regular}. 

Any regular boundary fixed point of a semigroup, other than the Denjoy--Wolff point, is characterized as \textit{repelling}. If $\sigma$ is a repelling fixed point of $(\phi_t)$, then there exists a $\lambda\in(-\infty,0)$ such that $\phi_t^{'}(\sigma)=e^{-\lambda t}$, for all $t\ge0$. This negative real number $\lambda$ is called the \textit{repelling spectral value} of $(\phi_t)$ at $\sigma$. 
Any non-regular boundary fixed point is characterized as \textit{super-repelling}.
\end{definition}

Following the notation of \cite[Chapter 13]{Booksem}, a \textit{backward orbit} of a semigroup $(\phi_t)$ is a continuous curve $\gamma: [0,+ \infty) \to \D$ that satisfies $\phi_s( \gamma(t))= \gamma(t-s)$, for all $t\ge s\ge0$.  The point $\gamma(0)$ is called the \textit{starting point} of $\gamma$.
A backward orbit is said to be \textit{regular} if $$\limsup_{t \to +\infty} d_{\D} (\gamma(t), \gamma(t+1)) < + \infty. $$
If the above limit is infinite, then we say that $\gamma$ is \textit{non-regular}.
We distinguish the following cases concerning the asymptotic behavior of regular backward orbits; the reader may refer to \cite[Chapter 13]{Booksem} for the detailed theory of backward orbits.

If $(\phi_t)$ is an elliptic semigroup then a backward orbit is either identical to the Denjoy--Wolff point $\tau\in\mathbb{D}$ or converges to a boundary fixed point of the semigroup. In the latter case, if the backward orbit is also regular, then it converges non-tangentially to a repelling fixed point.
If $(\phi_t)$ is a non-elliptic semigroup then a backward orbit necessarily converges to a boundary fixed point. If, in addition, the backward orbit is regular, then it converges either tangentially to the Denjoy--Wolff point or non-tangentially to a repelling fixed point. 

The collection of all backward orbits along with the forward trajectories of their starting points form the so-called \textit{backward invariant set} $\mathcal{W}$ of $(\phi_t)$. 

\begin{definition}\cite[Definition 13.4.1]{Booksem} 
\normalfont
Let $(\phi_t)$ be a one-parameter semigroup in $\mathbb{D}$.
A non-empty connected component $\Delta$ of the interior of $\mathcal{W}$ is called a \textit{petal}. 
\end{definition}

The restriction of every $\phi_t, t\ge0,$ in a petal $\Delta$ is an automorphism (i.e. $\phi_t(\Delta)=\Delta)$. Moreover, $\tau \in \partial \Delta$ and for every $z \in \Delta $ the curve $[0,+ \infty) \ni t  \mapsto \phi_{t}^{-1}(z)$ is a regular backward orbit for $(\phi_t)$. Without loss of generality, by denoting $\phi_{-t}(z):= \phi_t^{-1}(z)$, for all $z\in \Delta$, we can extend the semigroup in such a way that $t \in \mathbb{R}$. Hence, for the rest of this article, backward orbits will be defined for $t\le0$.

This extension is further applied at the images through the Koenigs function $h$. It can be easily seen that given $z\in\mathbb{D}$, in the case of non-elliptic semigroups, the image through $h$ of the backward orbit with starting point $z$ is actually a horizontal half-line emanating from $h(z)$ towards $\infty$ in the negative direction (i.e. with constant imaginary part and decreasing real part). In the case of elliptic semigroups, such a backward orbit is mapped through $h$ to a half-spiral reaching towards $\infty$.

For a semigroup $(\phi_t)$ with a petal $\Delta$, every point of $\Delta$ is actually the starting point of some regular backward orbit. Conversely, no point outside of a petal can be the starting point of a regular backward orbit. Therefore, the study of regular backward orbits is closely related to that of petals. In order to render this study simpler and more efficient, we once again turn to the Koenigs function $h$ of $(\phi_t)$. 

The image of a petal $\Delta$ under $h$ is always a maximal domain in $\Omega$. More specifically, for a non-elliptic semigroup, $h(\Delta)$ is
\begin{enumerate}
    \item[(i)] either a maximal horizontal strip, in which case the petal is characterized as \textit{hyperbolic},
    \item[(ii)] or a maximal horizontal half-plane, in which case the petal is characterized as \textit{parabolic}.
\end{enumerate}
Considering the one-to-one correspondence through $h$ between points of $\partial\mathbb{D}$ and prime ends of $\Omega$, it is deduced that all backward orbits contained inside a hyperbolic petal converge to the same repelling fixed point, while those inside a parabolic petal converge to the Denjoy--Wolff point. Moreover, in the former case, the width of this maximal horizontal strip depends on the repelling spectral value $\lambda$ of the repelling fixed point. In particular, the width is equal to $-\frac{\pi}{\lambda}$.

In the case where $(\phi_t)$ is an elliptic semigroup with Denjoy--Wolff point $\tau \in \D$ and spectral value $\mu$, the
image $h(\Delta)$ of a petal $\Delta$ is a maximal $\mu$-spirallike sector in $\Omega$ of center $e^{i \theta_0}$, for some $\theta_0 \in [-\pi,\pi)$, and amplitude $2a:=-\frac{|\mu|^2 \pi}{\lambda \RE \mu}$, where $\lambda$ is the repelling spectral value of the repelling fixed point to which every backward orbit contained in $\Delta$ converges. More concretely, 
$$h(\Delta) = \Spir \left[\mu , 2a, \theta_0 \right]: = \{e^{t \mu+i \theta} : t \in \mathbb{R}, \theta \in (-\alpha+\theta_0, \alpha+\theta_0)\}. $$
Every petal of an elliptic semigroup is characterized as \textit{hyperbolic}. The similarity between the two types of hyperbolic petals we described, relies on the fact that both of them contain a repelling fixed point in their boundaries.

\begin{lemma}\label{lem:mapofpetal}
    Suppose $(\phi_t)$ is an elliptic semigroup of holomorphic self maps of the unit disk $\D$ with spectral value $\mu \in \mathbb{C} $, $\RE \mu>0$. Denote by $h$ the associated Koenigs function of $(\phi_t)$. Let $\Delta$ be a petal of $(\phi_t)$ corresponding to the repelling fixed point $\sigma$ with repelling spectral value $\lambda$. Consider the function $g(z):= \frac{e^{-i \Arg \mu}}{\cos(\Arg \mu)} \Log z+ i (\alpha -\theta_0)$, where $a:=-\frac{|\mu|^2 \pi}{2 \lambda \RE \mu}$, $\theta_0 \in [-\pi,\pi)$, and $e^{i \theta_0}$ is the center of $h(\Delta)$. Then $g$ maps $h(\Delta)$ onto the horizontal strip $\{z\in \mathbb{C}: 0<\IM z<2 \alpha\}$. 
\begin{proof}
First note that $\cos(\Arg \mu)>0$, since $\RE \mu>0$.  Elementary calculations lead to $\Log e^{t \mu+i \theta} = t \mu+ i \theta$ and thus, spirals of $h(\Delta)$ are mapped onto lines parallel to the line joining $0$ with $\mu$.
Hence $$ g(h(\Delta)) = \bigcup_{\theta \in (-\alpha + \theta_0 , \alpha +\theta_0) } \left\{ t \frac{|\mu|}{\cos(\Arg \mu)}   + \tan(\Arg \mu) \theta + i (\theta+ \alpha -\theta_0): t \in \mathbb{R} \right\} $$
and since $\theta +\alpha -\theta_0 \in (0,2 \alpha) $, then $g$ maps $h(\Delta)$ onto the horizontal strip ${\{z \in \mathbb{C} : 0< \IM z < 2 \alpha \} }$.    
\end{proof}

\end{lemma}

Last but not least, independent of the type of the one-parameter semigroup, there exists a one-to-one correspondence between hyperbolic petals and repelling fixed points. We say that the repelling fixed point $\sigma$ corresponds to the hyperbolic petal $\Delta$ of $(\phi_t)$ (or $\Delta$ corresponds to $\sigma$, respectively) if all backward orbits lying in $\Delta$ converge to $\sigma$.

Keeping in mind the asymptotic behavior of backward orbits and the shape of petals, we can now discuss non-regular backward orbits. A non-regular backward orbit for one-parameter semigroups can fall into one of the following three cases:
\begin{enumerate}[(i)]
    \item it is either part of the boundary of a hyperbolic petal, in which case it converges tangentially to a repelling fixed point of the semigroup,
    \item it is either part of the boundary of a parabolic petal in which case it converges tangentially to the Denjoy--Wolff point of the semigroup, a situation that can arise solely in parabolic semigroups, 
    \item or it converges to a super-repelling fixed point of the semigroup, in which case the convergence can be either tangential or non-tangential (or an amalgamation of the two if we consider subsequences).
    \end{enumerate}

\subsection{ Infinitesimal Generators}

One-parameter semigroups are closely related with dynamical systems. In fact, for every $(\phi_t)$, there exists a unique holomorphic function $G: \D \to \mathbb{C}$, such that
$$\frac{\partial\phi_t(z)}{\partial t}=G(\phi_t(z)), \quad z\in\mathbb{D}, \quad t\ge0.$$
The function $G$ is called the \textit{infinitesimal generator} of the semigroup. 
In the course of the proofs, we exclusively need the following lemma regarding infinitesimal generators. It is a combination of results from Theorem 12.2.5 and Corollary 12.2.6 in \cite{Booksem}.

\begin{lemma}\label{lem:generator}
Let $(\phi_t)$ be a semigroup in $\D$ with associated infinitesimal generator $G$. Suppose $\sigma \in \partial \D$ is a repelling fixed point of $(\phi_t)$ with repelling spectral value $\lambda\in(-\infty,0)$. Then:
\begin{enumerate}
    \item[\textup{(i)}]  $$\RE \left( \frac{ \sigma G(z)}{(\sigma -z)^2} \right) \geq \frac{\lambda}{2} \cdot \frac{1-|z|^2}{|\sigma-z|^2}, \quad \text{for all } z \in \D;$$ 
    \item[\textup{(ii)}] $$-\lambda = \angle \lim_{z\to \sigma} \frac{G(z)}{z-\sigma}; $$ 
    \item[\textup{(iii)}] there exists a unique holomorphic mapping $p:\D \to \overline{\mathbb{H}}$ with $$\angle \lim_{z\to \sigma} (z-\sigma)p(z) =0, $$
    such that
    $$G(z ) =(\overline{\sigma}z-1)(z-\sigma) \left[p(z) +\frac{\lambda}{2} \cdot \frac{\sigma + z}{\sigma-z} \right], \quad\text{ for all } z\in\mathbb{D}.$$
\end{enumerate}
\end{lemma}

Further information on infinitesimal generators of one-parameter semigroups and their connection to the associated Koenigs functions can be found in \cite[Chapter 10]{Booksem} as well as in references therein.

\section{Speeds of Convergence for Petals}\label{sec:speeds}
 
Before we proceed to the definitions of speeds of convergence along backward orbits, we need some preliminary propositions. Recall that for every $z\in\Delta$, where $\Delta$ is a petal of a semigroup $(\phi_t)$, the curve with $(-\infty,0]\ni t\mapsto\phi_t(z)=\phi_{-t}^{-1}(z)$ is a backward orbit. 

\begin{proposition}\label{prop:total}
Let $(\phi_t)$ be a semigroup in $\D$ and $\Delta$ be a petal of $(\phi_t)$. Let $z,w\in\Delta$. Then
$$|d_\mathbb{D}(z,\phi_t(z))-d_\mathbb{D}(w,\phi_t(w))|\le2d_\Delta(z,w), \quad t\le0.$$
\begin{proof}
First, suppose that $(\phi_t)$ is non-elliptic and let $h$ be its Koenigs function. Set $\Omega=h(\mathbb{D})$. Using consecutively the conformal invariance of the hyperbolic distance and the triangle inequality, we have for $t\le0$,
\begin{eqnarray*}
d_\mathbb{D}(z,\phi_t(z))-d_\mathbb{D}(w,\phi_t(w))&=&d_\Omega(h(z),h(z)+t)-
d_\Omega(h(w),h(w)+t) \\
&\le&d_\Omega(h(z),h(w))+d_\Omega(h(z)+t,h(w)+t).
\end{eqnarray*}
However, $h(\Delta)\subset\Omega-t\subset\Omega$, for all $t\le0$. Therefore, continuing with the previous inequalities, we get
\begin{eqnarray*}
d_\Omega(h(z),h(w))+d_\Omega(h(z)+t,h(w)+t)&\le&d_{h(\Delta)}(h(z),h(w))+d_{\Omega-t}(h(z),h(w))\\
&\le&2d_{h(\Delta)}(h(z),h(w))\\
&=&2d_\Delta(z,w).
\end{eqnarray*}
Interchanging the roles of $z$ and $w$ we get a similar inequality, which in turn leads to the desired result. The proof for elliptic semigroups is almost identical, albeit with some minor modifications because of the different form of the Koenigs function.
\end{proof}
\end{proposition}
\begin{remark}
    It is easy to check that the result of this proposition does not remain true if $z$ and $w$ lie on different petals or if any of the two lies outside of any petal. In such cases, the above difference between what will eventually be the total speeds, does not remain bounded. 
\end{remark}

In the study of speeds in the forward dynamics, the role of the geodesic was assumed by the diameter with one end on the Denjoy--Wolff point $\tau$ of the non-elliptic semigroup, because every trajectory converges to this point. The choice of the diameter, instead of any other geodesic landing at $\tau$, is justified as $0$ is the starting point of the trajectory used in ``forward'' speeds and the diameter is the only such geodesic that passes through $0$. In this way, $v(0)=v^o(0)=v^T(0)=0$.
However, this is not the case when it comes to backward dynamics, since the selection of the geodesic fluctuates depending on the petal. Once again, let $(\phi_t)$ be a semigroup with a petal $\Delta$ and let $z\in\Delta$. We denote by $\eta:(-1,1)\to\mathbb{D}$ the geodesic for the hyperbolic distance $d_\mathbb{D}$ such that:
\begin{enumerate}
    \item[(i)] if $\Delta$ is hyperbolic, then $\eta(0)=z$ and $\lim\limits_{r\to1^-}\eta(r)=\sigma$, where $\sigma\in\partial\mathbb{D}$ is the repelling fixed point of $(\phi_t)$ where all backward orbits contained in $\Delta$ converge;
    \item[(ii)] if $\Delta$ is parabolic, then $\eta(0)=z$ and $\lim\limits_{r\to1^-}\eta(r)=\tau$, where $\tau\in\partial\mathbb{D}$ is the Denjoy--Wolff point of $(\phi_t)$.
\end{enumerate}

\begin{proposition}\label{prop:ort&tang}
Let $(\phi_t)$ be a semigroup in $\D$ and $\Delta$ be a petal of $(\phi_t)$. Let $z,w\in\Delta$ and suppose that $\eta:(-1,1)\to\mathbb{D}$ is the geodesic described above, with $\eta(0)=z$. Then
$$|d_\mathbb{D}(z,\pi_\eta(\phi_t(z)))-d_\mathbb{D}(w,\pi_\eta(\phi_t(w)))|\le2d_\Delta(z,w), \quad t\le0$$
and 
$$|d_\mathbb{D}(\phi_t(z),\eta)-d_\mathbb{D}(\phi_t(w),\eta)|\le 2d_\Delta(z,w), \quad t\le0.$$
\begin{proof}
For the sake of brevity, for $t\le0$ we write
$$z_t=\pi_\eta(\phi_t(z)) \quad \text{and} \quad w_t=\pi_\eta(\phi_t(w)).$$
The proof is almost identical with that of Proposition \ref{prop:total}, but there is a key argument that renders the presentation of this proof necessary. We have
\begin{eqnarray*}
d_\mathbb{D}(z,z_t)-d_\mathbb{D}(w,w_t)&=&d_\Omega(h(z),h(z_t))-d_\Omega(h(w),h(w_t))\\
&\le&d_\Omega(h(z),h(w))+d_\Omega(h(z_t),h(w_t))\\
&\le&d_{h(\Delta)}(h(z),h(w))+d_\mathbb{D}(z_t,w_t).
\end{eqnarray*}
Using \cite[Proposition 3.3]{brspeeds}, we get
$$d_\mathbb{D}(z_t,w_t)\le d_\mathbb{D}(\phi_t(z),\phi_t(w)), \quad \text{for all }t\le0.$$
Combining all the above inequalities and continuing exactly as in the proof of the Proposition \ref{prop:total}, we deduce the first desired result.
The proof for the second desired inequality follows the same steps so we skip it in order to avoid repetitiveness. 
\end{proof}
\end{proposition}

Once more, the bounds established demonstrate that the asymptotic behavior of all the aforementioned quantities does not depend on the choice of the starting point provided we stay inside a petal. All the propositions that we proved point to the fact that the key difference between the two kinds of speeds is the following: while the speed of a forward trajectory can be globally defined as the speed of the semigroup, this is not true about the speed along a backward orbit. On the contrary, the speed along a regular backward orbit can be locally defined as the speed of the respective petal.

\begin{definition}\label{def:speeds}
    \normalfont
    Let $(\phi_t)$ be a semigroup in $\mathbb{D}$ and let $\Delta$ be a petal of $(\phi_t)$. Let $z\in\Delta$ and suppose that $\eta:(-1,1)\to\mathbb{D}$ is the geodesic described earlier, depending on $z$ and the type of the petal. We call \textit{total speed} of $\Delta$ the function
    $$v_\Delta(t)=d_\mathbb{D}(z,\phi_t(z)), \quad t\le0.$$
    We call \textit{orthogonal speed} of $\Delta$ the function
    $$v^o_\Delta(t)=d_\mathbb{D}(z,\pi_\eta(\phi_t(z))), \quad t\le0.$$
    We call \textit{tangential speed} of $\Delta$ the function
    $$v^T_\Delta(t)=d_\mathbb{D}(\phi_t(z),\eta)=d_\mathbb{D}(\phi_t(z),\pi_\eta(\phi_t(z))), \quad t\le0.$$ 
\end{definition}

A first useful result that correlates the three speeds stems directly from Bracci's Pythagoras' Theorem \cite[Proposition 3.4]{Booksem}.

\begin{corollary}
Let $(\phi_t)$ be a semigroup in $\D$ and $\Delta$ be a petal of $(\phi_t)$. Then
$$v^o_\Delta(t)+v^T_\Delta(t)-\frac{1}{2}\log2\le v_\Delta(t)\le v^o_\Delta(t)+v^T_\Delta(t),$$
for all $t\le0$.
\end{corollary}

\section{Total Speed of Petals}\label{sec:totalspeed}

In this section, we examine the convergence of the total speed of a petal $\Delta$. The asymptotic behavior of the speeds of a petal depends greatly on the type of the petal.

We first deal with hyperbolic petals which may exist in both elliptic and non-elliptic semigroups. Suppose $(\phi_t)$ is a one-parameter semigroup. Let $\sigma$ be a repelling fixed point of $(\phi_t)$ with spectral value $\lambda <0$. Denote by $\Delta$ the associated hyperbolic petal. Our goal is to prove results, analogue to that of relation \eqref{speeds2}. As a matter of fact, we obtain the following result. 

\begin{theorem}\label{thm:divratehyperbolic}
Let $(\phi_t)$ be a semigroup in $\D$. Suppose $\Delta$ is a hyperbolic petal of $(\phi_t)$ corresponding to the repelling fixed point $\sigma$. Then
\begin{equation}
    \lim_{t \to -\infty} \frac{v_{\Delta}(t)}{t}= \frac{\lambda}{2}, 
\end{equation}
where $\lambda $ is the repelling spectral value of $\sigma$.
\end{theorem}

For the proof of Theorem \ref{thm:divratehyperbolic}, we need information on the asymptotic behavior of the hyperbolic distance inside the hyperbolic petal $\Delta$. Restating \cite[Lemma 13.5.1]{Booksem} with a slight re-parametrization, the following lemma arises.

\begin{lemma}\label{lem:hypdistpetal}
Let $(\phi_t)$ be a semigroup in $\D$. Suppose $\Delta$ is a hyperbolic petal of $(\phi_t)$ corresponding to the repelling fixed point $\sigma$. Then  
\begin{equation}
    \lim_{t \to -\infty} \frac{d_{\Delta}(z, \phi_t(z))}{t}= \frac{\lambda}{2}, \quad \text{for all } z \in \Delta,  
\end{equation}
where $\lambda $ is the repelling spectral value of $(\phi_t)$ at $\sigma$.
\end{lemma}

Theorem \ref{thm:divratehyperbolic} suggests that the same limit is true when we use the hyperbolic distance in the unit disk instead of restricting to the hyperbolic geometry of the petal. The main idea for the proof is that while the backward orbit converges to the repelling fixed point, the hyperbolic geometry of the unit disk becomes similar to that of the petal.

\begin{proof}[\bf Proof of Theorem \ref{thm:divratehyperbolic}]
Let $z\in \Delta$. Due to domain monotonicity, we can easily observe that 
\begin{equation}\label{eq:proof_lowerbound}
    \liminf_{t \to -\infty} \frac{d_{\D} (z, \phi_t(z))}{t} \geq \lim_{t \to -\infty} \frac{d_{\Delta}(z, \phi_t(z)) }{t}= \frac{\lambda }{2}. 
\end{equation}
For the other way inequality, we use l'H\^{o}pital's Rule in order to estimate the limit. From the generalized l'H\^{o}pital's Rule, we obtain 
\begin{equation}\label{eq:lhospital}
    \limsup_{t \to -\infty} \frac{d_{\D}(z, \phi_t(z))}{t} \leq \limsup_{t \to -\infty}\frac{\partial}{\partial t} d_{\D}(z, \phi_t(z))  .
\end{equation}
Let us denote by $d_z(t)$ the hyperbolic distance $d_{\D}(z, \phi_t(z))$ and by $\rho_z (t)$ the pseudo-hyperbolic distance $\rho_{\D}(z, \phi_t(z))$. 
We can calculate the derivative of the hyperbolic distance w.r.t. $t$ using its representation in terms of the pseudo-hyperbolic distance; $2d_z (t) =\log (1+\rho_z (t)) - \log(1-\rho_z (t))$.  With elementary calculations, we obtain
\begin{equation}\label{eq:derivativeofpseudo}
2 \frac{\partial}{\partial t} d_z (t) =\frac{1}{\rho_z (t) (1-\rho_z (t)^2)} \frac{\partial }{\partial t} \rho_z (t)^2. 
\end{equation}
The pseudo-hyperbolic distance can be explicitly expressed in the following way $$\rho_z (t)^2 = \left| \frac{z-\phi_t(z)}{1-\bar{z} \phi_t(z)} \right|^2 = \frac{|z|^2 +|\phi_t(z)|^2 - 2 \RE(\bar{z} \phi_t(z)) }{1+|z|^2 \cdot |\phi_t(z)|^2 -2 \RE(\bar{z} \phi_t(z)) }. $$
In order to calculate the derivative in \eqref{eq:derivativeofpseudo}, we use the infinitesimal generator $G$ of $(\phi_t)$, since $\frac{\partial\phi_t(w)}{\partial t}=G(\phi_t(w))$, for all $w\in\Delta$ and $t\le0$. Hence we are led to 
\begin{equation}\label{eq:???}
    \rho_z (t) \frac{\partial}{\partial t} d_z (t)  = \frac{ (1+|z|^2-2\RE(\bar{z} \phi_t(z)))\RE(G(\phi_t(z)) \overline{\phi_t(z)} )}{{ (1-|\phi_t(z)|^2) |1-\bar{z} \phi_t(z)|^2}} -\frac{\RE(G(\phi_t(z)) \bar{z})}{ |1-\bar{z} \phi_t(z)|^2} .
\end{equation}
We note with elementary calculations that $$  \liminf_{t\to -\infty} \frac{\RE(G(\phi_t(z)) \bar{z})}{ |1-\bar{z} \phi_t(z)|^2} =0, $$
since $G(\sigma)=0$. We write $$ \RE(G(\phi_t(z)) \overline{\phi_t(z)} ) = \RE \left( \frac{G(\phi_t(z))\sigma}{(\sigma-\phi_t(z))^2} \overline{\sigma} \overline{\phi_t(z)} (\sigma-\phi_t(z))^2 \right)$$ and further observe that 
$$ \overline{\sigma} \overline{\phi_t(z)} \frac{(\sigma-\phi_t(z))^2}{1-|\phi_t(z)|^2} = 1-\overline{\sigma}\phi_t(z) - \frac{|\sigma-\phi_t(z)|^2}{1-|\phi_t(z)|^2}. $$
As a result, 
$$ \frac{\RE(G(\phi_t(z)) \overline{\phi_t(z)} )}{1-|\phi_t(z)|^2} = \RE \left( \frac{G(\phi_t(z))}{\sigma - \phi_t(z)}\right)  - \frac{|\sigma - \phi_t(z)|^2}{1-|\phi_t(z)|^2} \RE \left( \frac{G(\phi_t(z)) \sigma}{(\sigma - \phi_t(z))^2} \right) $$
From Lemma \ref{lem:generator} we obtain 
\begin{equation}\label{eq:proof_infgener1}
\RE \left( \frac{G(\phi_t(z))}{\sigma - \phi_t(z)} \right)  =  - \RE \left( (\overline{\sigma} \phi_t(z)-1) p(\phi_t(z)) \right) + \frac{\lambda}{2} \RE( \overline{\sigma} (\sigma+\phi_t(z)) ),    
\end{equation}
where $p$ is a holomorphic mapping of positive real part. 
Moreover, again from Lemma \ref{lem:generator}
\begin{equation}\label{eq:proof_infgener2}
    - \frac{|\sigma - \phi_t(z)|^2}{1-|\phi_t(z)|^2} \RE \left( \frac{G(\phi_t(z)) \sigma}{(\sigma - \phi_t(z))^2} \right) \leq -\frac{\lambda}{2} \frac{|\sigma - \phi_t(z)|^2}{1-|\phi_t(z)|^2} \frac{1-|\phi_t(z)|^2}{|\sigma - \phi_t(z)|^2} = - \frac{\lambda}{2}
\end{equation}
Hence returning back to \eqref{eq:lhospital} and combining \eqref{eq:proof_infgener1} and \eqref{eq:proof_infgener2}, it follows that 
\begin{eqnarray*}
\limsup_{t\to-\infty} \frac{d_z (t)}{t} &\leq & \limsup_{t \to -\infty} \frac{1+|z|^2-2 \RE(\bar{z} \phi_t(z))}{|1-\bar{z} \phi_t(z)|^2} \left[ - \RE \left( (\overline{\sigma} \phi_t(z)-1) p(\phi_t(z)) \right) \right.\nonumber\\
&\quad& \quad \quad \left. + \frac{\lambda}{2} \RE( \overline{\sigma} (\sigma+\phi_t(z)) ) -\frac{\lambda}{2} \right] \nonumber \\
&= & \lambda - \frac{\lambda}{2} =\frac{\lambda}{2}.
\end{eqnarray*}
Taking also \eqref{eq:proof_lowerbound} into account, the desired limit occurs. 
\end{proof}

Let us recall that parabolic petals exist only in the case of parabolic one-parameter semigroups. Suppose that $(\phi_t)$ is a parabolic semigroup of holomorphic self-maps of $\D$ with Denjoy--Wolff point $\tau \in \partial\mathbb{D}$. Further suppose that $h$ is the associated Koenigs function of the semigroup. Every backward orbit contained inside a parabolic petal $\Delta$ of the semigroup converges to $\tau$. We obtain the following result on the asymptotic behavior of the total speed of a parabolic petal. 

\begin{theorem}\label{thm:bounds}
Let $(\phi_t)$ be a semigroup in $\mathbb{D}$ with a parabolic petal $\Delta$. Then, there exists an absolute constant $c\ge1$ such that
     $$\frac{1}{c}\cdot\log|t|\le v_{\Delta}(t)\le c\cdot\log|t|,$$
     for all $t<-1$.
\begin{proof}
Since $\Delta$ is parabolic, $(\phi_t)$ is parabolic as well and thus $\Omega=h(\mathbb{D})$ is convex in the positive direction, while $h(\Delta)$ is a maximal horizontal half-plane inside $\Omega$. Let $z\in\Delta$. By the conformal invariance and the monotonicity property of hyperbolic distance, for $t\le0$ we have
$$v_\Delta(t)=d_\mathbb{D}(z,\phi_t(z))=d_\Omega(h(z),h(z)+t)\le d_{h(\Delta)}(h(z),h(z)+t).$$ 
Without loss of generality, we assume that $h(\Delta)$ is the upper half-plane and $h(z)=i$ (in a different case, this can be achieved by a composition with a translation). Through a rotation, we turn to the right half-plane $\mathbb{H}$ and with the use of \eqref{eq:hypdiststripandhalfplane}, we obtain
   $$             d_\mathbb{D}(z,\phi_t(z)) \le d_{h(\Delta)}(i,i+t)=d_\mathbb{H}(1,1-it) = \arctanh \frac{-t}{\sqrt{4+t^2}}.$$
Therefore, for $t<-1$, we get
    \begin{equation}
        \limsup\limits_{t\to-\infty}\frac{d_\mathbb{D}(z,\phi_t(z))}{\log|t|}\le  \lim\limits_{t\to-\infty}\frac{\arctanh \frac{-t}{\sqrt{4+t^2}}}{\log(-t)}.
     \end{equation}
     Using l'H\^{o}pital's Rule, we can compute the last limit and find that
     \begin{equation}
        \limsup\limits_{t\to-\infty}\frac{v_\Delta(t)}{\log|t|}\le  1,
     \end{equation}
     where we have avoided the explicit mention of some redundant calculations.
        We now search for a reverse inequality. From Lemma \ref{lem:distancelemma}, it follows that
        $$d_\mathbb{D}(z,\phi_t(z))=d_\Omega(h(z),h(z)+t)\ge\frac{1}{4}\log\left(1+\frac{|h(z)-(h(z)+t)|}{\min\{\delta_\Omega(h(z)),\delta_\Omega(h(z)+t)\}}\right),$$
        for all $t\le0$. Since the associated planar domain $\Omega$ is convex in the positive direction, the boundary distance $\delta_\Omega(h(z)+t)$ is an non-decreasing function of $t\le0$. As a result,
        $$d_\mathbb{D}(z,\phi_t(z))\ge\frac{1}{4}\log\left(1+\frac{|t|}{\delta_\Omega(h(z)+t)}\right),$$ for all $t\le0$. 
        However, the set $\{h(z)+t:t\le0\}$ is the trace of a backward orbit through the Koenigs function $h$ and thus, it is contained in a maximal horizontal half-plane. This means that the value of $\delta_\Omega(h(z)+t)$ remains bounded and strictly positive, as $t\to-\infty$. We can write  $\lim\limits_{t\to-\infty}\delta_\Omega(h(z)+t)=:d\in(0,+\infty)$. Let $\epsilon>0$. Then, there exists $t_0\le-1$ such that $\delta_\Omega(h(z)+t)\le d+\epsilon$, for all $t\le t_0$. Consequently and with a usage of l'H\^{o}pital's Rule, we find
        \begin{eqnarray*}
        \liminf\limits_{t\to-\infty}\frac{d_\mathbb{D}(z,\phi_t(z))}{\log|t|}&\ge&\frac{1}{4}\liminf\limits_{t\to-\infty}\frac{\log\left(1+\frac{-t}{\delta_\Omega(h(z)+t)}\right)}{\log(-t)}\\
        &\ge&\frac{1}{4}\liminf\limits_{t\to-\infty}\frac{\log\left(1-\frac{t}{d+\epsilon}\right)}{\log(-t)} = \frac{1}{4}.
        \end{eqnarray*}
        All in all, we get
        $$\frac{1}{4}\le\liminf\limits_{t\to-\infty}\frac{v_\Delta(t)}{\log|t|}\le\limsup\limits_{t\to-\infty}\frac{v_\Delta(t)}{\log|t|}\le1.$$
        This last relation certifies the existence of the required $c\ge1$ and we get the desired result.
     \end{proof}
     \end{theorem}

    The statement of Theorem \ref{thm:bounds} guarantees, in other words, that the total speed of a parabolic petal is actually comparable to $\log|t|$.
    
    \begin{corollary}\label{cor:distinparabpetal}
    Let $(\phi_t)$ be a semigroup in $\mathbb{D}$ and $\Delta$ a parabolic petal of $(\phi_t)$. Then,
    $$\lim\limits_{t\to-\infty}\frac{v_\Delta(t)}{t}=0.$$
    \end{corollary}
    \begin{remark}
    Remember that if $\Delta$ is a hyperbolic petal of a non-elliptic semigroup with associated repelling fixed point $\sigma$, then the image $h(\Delta)$ is a maximal horizontal strip of width $-\frac{\pi}{\lambda}$, where $\lambda<0$ is the repelling spectral value of $\sigma$. When $\Delta$ is parabolic, the image is a maximal horizontal half-plane. In a way, we can say that $h(\Delta)$ is a horizontal strip of infinite width. Thus, the ``repelling'' spectral value can be said to be $0$, something that agrees with the actual spectral value of the semigroup at the Denjoy-Wolff point. Taking all the above into account, we may state that Corollary \ref{cor:distinparabpetal} is a natural extension of Theorem \ref{thm:divratehyperbolic}.
    \end{remark}

\section{Tangential and Orthogonal Speeds of Petals}\label{sec:orthogonal&tandspeed}

The asymptotic behavior of the tangential and orthogonal speeds of a petal is the main subject of this section. As outlined in the case of the total speed, there are two separate cases depending on the type of the petal. We obtain the following result concerning the asymptotic behavior of tangential speeds of petals. 

\begin{theorem}\label{thm:asympttang}
Let $(\phi_t)$ be a semigroup in $\mathbb{D}$. Suppose that $\Delta$ is a petal of $(\phi_t)$. The following are true. 
\begin{enumerate}
    \item[\textup{(i)}] If $\Delta$ is hyperbolic, then $\limsup\limits_{t\to-\infty}v^T_\Delta(t)<+\infty.$
    \item[\textup{(ii)}] If $\Delta$ is parabolic, then $\lim\limits_{t\to-\infty}v^T_\Delta(t)=+\infty.$
\end{enumerate}
\begin{proof}
The proof is long and mainly relies on estimates of harmonic measure. We divide it in three distinct cases:

\textbf{Case 1} \textit{ $\Delta$ is hyperbolic and $(\phi_t)$ is non-elliptic} : 
Let $\tau$ be the Denjoy--Wolff point of the semigroup, $h$ its Koenigs function and $\Omega=h(\mathbb{D})$. Let $z\in\Delta$. The hyperbolic petal $\Delta$ corresponds to a repelling fixed point $\sigma\in\partial\mathbb{D}$ of the semigroup and so, the backward orbit emanating from $z$ converges to $\sigma$. Suppose that $\eta:(-1,1)\to\mathbb{D}$ is the geodesic of the hyperbolic distance in the unit disk with $\eta(0)=z$ and $\lim\limits_{r\to1^-}\eta(r)=\sigma$. Then
$$v^T_\Delta(t) = d_\mathbb{D}(\phi_t(z),\eta) = d_\Omega(h(z)+t,h\circ\eta) =d_\Omega(h(z)+t,x_t),$$
for all $t\le0$, where $x_t\in h(\eta((-1,1)))$ is the point such that
$$d_\Omega(h(z)+t,x_t)=d_\Omega(h(z)+t,h\circ\eta)=\inf\limits_{r\in(-1,1)}d_\Omega(h(z)+t,h(\eta(r))).$$
The image $h(\Delta)$ is a maximal horizontal strip $S$ in $\Omega$ which contains $h(z)+t$, for all $t\le0$. Our first aim is to prove that the points $x_t$ on the geodesic $h\circ\eta$ of $\Omega$ are eventually contained inside $S$ as well.

Suppose that there exists a strictly decreasing sequence $\{t_n\}\subset(-\infty,0]$ with $t_n \xrightarrow{n\to+\infty}-\infty$ such that $x_{t_n}\in\Omega\setminus S$, for every $n\in\mathbb{N}$. Without loss of generality, we can assume by potentially extracting a subsequence that $\IM x_{t_n}> \IM h(z)$, for every $n\in\mathbb{N}$. The fact that the horizontal strip $S$ is maximal in $\Omega$ leads to 
$$\lim\limits_{n\to+\infty}\delta^+_\Omega(x_{t_n})=0\quad \text{and} \quad \lim\limits_{n\to+\infty}\delta^-_\Omega(x_{t_n})>0,$$
where $\delta^+_\Omega(z):=\dist(z,\partial\Omega^+)$ for $\partial\Omega^+=\{\zeta\in\partial\Omega: \IM \zeta> \IM h(z)\}$ and $\delta^-_\Omega(z):=\dist(z,\partial\Omega^-)$ for $\partial\Omega^-=\{\zeta\in\partial\Omega:\IM\zeta< \IM h(z)\}$.
Trivially then,
\begin{equation}\label{ratio}
\lim\limits_{n\to+\infty}\frac{\delta^+_\Omega(x_{t_n})}{\delta^-_\Omega(x_{t_n})}=0.
\end{equation}
Extracting again a subsequence, we may assume that $\delta_\Omega^+(x_{t_n})<\delta_\Omega^-(x_{t_n})$, for all $n\in\mathbb{N}$. Consider $p_n^+$ to be the point on $\partial\Omega^+$ such that $|x_{t_n}-p_n^+|=\delta_\Omega^+(x_{t_n})$ and $p_n^-$ the corresponding point on $\partial\Omega^-$. Set $\gamma_n^+=[x_{t_n},p_n^+)$ and $\gamma_n^-=[x_{t_n},p_n^-)$, where the symbols signify the respective line segments. Obviously, $\gamma_n^+,\gamma_n^-\subset\Omega$. Next, consider $D_n:=D(x_{t_n},\delta_\Omega^-(x_{t_n}))$ to be the disk of center $x_{t_n}$ and radius $\delta_\Omega^-(x_{t_n})>0$. By construction, $\gamma_n^+\subset D_n$, while $\gamma_n^-$ is a radius of $D_n$. Finally, denote by $\Omega_n$ the connected component of $\Omega\cap D_n$ containing $\gamma_n^+$ and set $E_n=\partial \Omega_n \setminus\partial\Omega^+$. Every curve that joins $\gamma_n^+$ to $\partial\Omega^-$ inside $\Omega$ definitely contains a subcurve joining $\gamma_n^+$ to $E_n$ inside $\Omega_n$. Turning to the definition of extremal distance and the first monotonicity property of extremal length, we get
$$\lambda(\gamma_n^+,\partial\Omega^-,\Omega)\ge\lambda(\gamma_n^+,E_n,\Omega_n),$$
for every $n\in\mathbb{N}$. Moreover, the family of curves joining $\gamma_n^+$ to $E_n$ inside $\Omega_n$ is evidently a subfamily of the family of curves joining $\gamma_n^+$ to $\partial(D_n)$ inside the disk $D_n$. Therefore, by Lemma \ref{extr-length-monoton}
$$\lambda(\gamma_n^+,E_n,\Omega_n)\ge\lambda(\gamma_n^+,\partial(D_n),D_n),$$
for every $n\in\mathbb{N}$. Using back-to-back a translation, a dilation and a rotation, we see that
$$\lambda(\gamma_n^+,\partial(D_n),D_n)=\lambda\left(\left[0,\frac{\delta_\Omega^+(x_{t_n})}{\delta_\Omega^-(x_{t_n})}\right),\partial\D,\D\right),$$
for every $n\in\mathbb{N}$. Taking into account \eqref{ratio}, we turn to the theory of Gr\"{o}tzsch rings (see \cite[Theorem 2.80]{ohtsuka}) to find out that
$$\lim\limits_{n\to+\infty}\lambda\left(\left[0,\frac{\delta_\Omega^+(x_{t_n})}{\delta_\Omega^-(x_{t_n})}\right),\partial\D,\D\right)=+\infty.$$
Tracing back our steps, it follows that
$$\lim\limits_{n\to+\infty}\lambda(\gamma_n^+,\partial\Omega^-,\Omega)=+\infty.$$
Finally, with the help of Beurling's estimate, we are able to reformulate this result in terms of harmonic measure and write
$$\lim\limits_{n\to+\infty}\omega(x_{t_n},\partial\Omega^-,\Omega)=0.$$
In view of Remark \ref{caratheodory}, the asymptotic behavior of the above quantity implies that the sequence $\{h^{-1}(x_{t_n})\}$ converges to $\sigma$ tangentially. However, $h\circ\eta$ is actually a geodesic for the hyperbolic distance in $\Omega$ and conformal mappings preserve angles. This signifies that $\{h^{-1}(x_{t_n})\}$ converges to $\sigma$ orthogonally. As a result, we are led to a contradiction! Consequently, there exists a $t_0\le0$ such that $x_t\in S$, for all $t\le t_0$. 

\begin{figure}[h]
    \centering
    \includegraphics[scale=0.45]{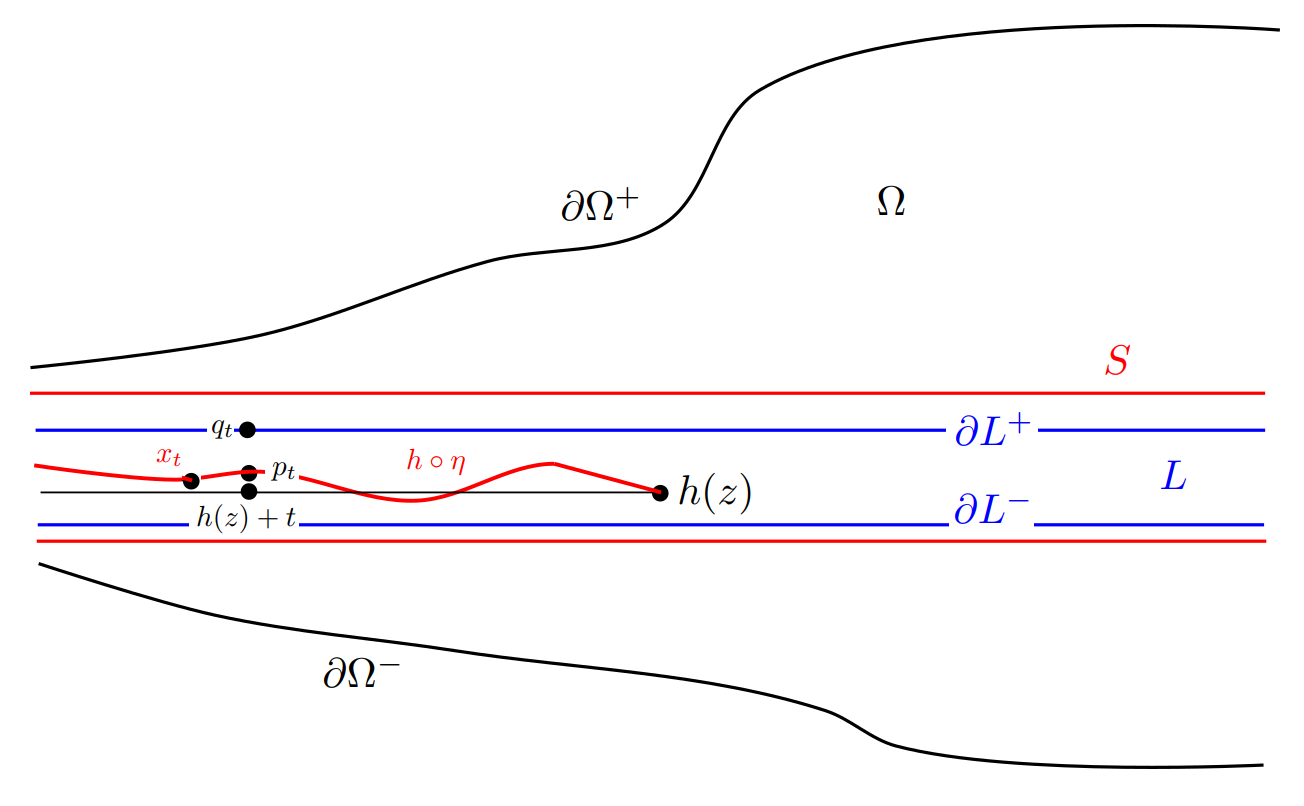}
    \caption{Sets of Theorem \ref{thm:asympttang}: \textbf{Case 1}}
    \label{hyperbolic-limsup}
\end{figure}

By the monotonicity property of hyperbolic distance, we get $v^T_\Delta(t)\le d_S(h(z)+t,x_t),$ for all $t\le t_0$.
 Insisting on our previous argument that $\{x_t\}$ is eventually contained inside $S$, we understand that there must exist another horizontal strip $L$ such that $L\subsetneq S$, $\partial L\cap\partial S=\emptyset$ and $\{x_t\}$ is eventually contained in this smaller strip $L$ (otherwise the limit at \eqref{ratio} would be $0$ or $+\infty$, leading once again to tangential convergence and a contradiction). Clearly, we can choose $L$ in such a way that $h(z)\in L$, too. For $t\le t_0$ consider $p_t$ the point on $h\circ\eta((-1,1))$ such that $ \RE p_t= \RE h(z)+t$. We know that there exists such a point for all $t\le t_0$, since $h(\eta(0))=h(z)$ and because $h\circ\eta$ converges to the prime end of $\Omega$ corresponding to $\sigma$ which is defined inside the strip $S$ and to the negative direction. In addition, the point $p_t$ is unique by Minda's Reflection Principle (see \cite{MindaReflection}). By definition
$$d_S(h(z)+t,x_t)\le d_S(h(z)+t,p_t).$$

Consider now $q_t$ to be the point of $\partial L^+$ if $ \IM p_t\ge \IM h(z)$ or the point of $\partial L^-$ if $ \IM p_t\le \IM h(z)$ such that $ \RE q_t= \RE h(z)+t$ (see Figure \ref{hyperbolic-limsup}). 

 However, the vertical arc joining $h(z)+t$ with $p_t$ and then with $q_t$ is a geodesic arc of the horizontal strip $S$. Due to the orientation, this implies that
$$v^T_\Delta(t)\le d_S(h(z)+t,p_t)\le d_S(h(z)+t,q_t),$$
for all $t\le t_0$. Nevertheless, horizontal translations leave the strip $S$ invariant and the points $h(z)+t$ and $q_t$ share the same real part. Therefore, the quantity $d_S(h(z)+t,q_t)$ can only take two distinct values $c^+$ and $c^-$, depending on if $q_t\in L^+$ or $q_t\in L^-$, respectively. This implies that
$$v^T_\Delta(t)\le\max\{c^+,c^-\},$$
for all $t\le t_0$ and thus, $\limsup\limits_{t\to-\infty}v^T_\Delta(t)<+\infty$.

\textbf{Case 2} \textit{$\Delta$ is hyperbolic and $(\phi_t)$ is elliptic} :  
Since $(\phi_t)$ is elliptic, $h(\Delta)$ is a maximal $\mu$-spirallike sector in $\Omega$, where $\mu$ is the spectral value of the semigroup. With the help of Lemma \ref{lem:mapofpetal}, $h(\Delta)$ can be mapped onto a strip. Working as in the previous case, we may prove the desired result. The only difference is that the images of the backward orbit and the geodesic through $g \circ h$ converge to $\infty$ in the positive direction.

\textbf{Case 3} \textit{$\Delta$ is parabolic} : 
The image $h(\Delta)$ is a maximal horizontal half-plane inside $\Omega$. Set $H=h(\Delta)$ and let $z\in\Delta$. Working similarly to the previous proof, we have
$$  v_\Delta^T(t) = d_\mathbb{D}(\phi_t(z),\eta) = d_\Omega(h(z)+t,h\circ\eta) = d_\Omega(h(z)+t,x_t),$$
for all $t\le0$, where by $x_t$ we denote again the point on $h\circ\eta$ that is hyperbolically closest to $h(z)+t$.
Utilizing Lemma \ref{lem:distancelemma}, we get
$$v_\Delta^T(t)\ge\frac{1}{4}\log\left(1+\frac{|h(z)+t-x_t|}{\min\{\delta_\Omega(h(z)+t),\delta_\Omega(x_t)\}}\right).$$ 
Our goal is to estimate every factor in the right hand side of the inequality. Due to the convexity in the positive direction of $\Omega$ and the maximality of $H$ in $\Omega$, we deduce that
$\lim\limits_{t\to-\infty}\delta_\Omega(h(z)+t)=c\in(0,+\infty)$.

We further prove that $\liminf\limits_{t\to-\infty}|h(z)+t-x_t|>0$. Pursuing a contradiction, let us assume that $\liminf\limits_{t\to-\infty}|h(z)+t-x_t|=0$. Then, there exists a strictly decreasing sequence $\{t_n\}\subset(-\infty,0]$ such that $\lim\limits_{n\to+\infty}t_n=-\infty$ and $\lim\limits_{n\to+\infty}|h(z)+t_n-x_{t_n}|=0$. In order for the last limit to be zero, it is necessary that $\lim\limits_{n\to+\infty} \RE x_{t_n}=-\infty$, while $ \IM x_{t_n}$ must remain bounded. Therefore, $\lim\limits_{n\to+\infty}\Arg x_{t_n}=\pi$. In addition, following the steps of the proof in \textbf{Case 1}, we can safely assume that $\{x_{t_n}\}$ is contained in $H$ as well.

For the remainder of the proof, we make use of the harmonic measure. Its conformal invariance allows us to assume, without loss of generality, that $H$ is the upper half-plane. Set $L^+=\{z\in\mathbb{C}: \Arg z=\pi\}\cup\{0\}$ and $L^-=\partial H\setminus L^+$. We first deal with the case when $\partial H\cap\partial\Omega=\emptyset$. Since $H\subset\Omega$, it must be true that $ \IM w<0$ for all $w\in\partial\Omega$. In particular, $\partial\Omega$ contains a connected component $E$ such that $\dist(E,L^+)=0$ because the two sets tend to intersect when approaching $\infty$ in the negative direction. Set $\zeta$ the point of $E$ closest to $0$ and $A$ the line segment joining $0$ and $\zeta$. Our construction, the maximality of $H$ in $\Omega$ and the convexity in the positive direction of $\Omega$ imply that $A\setminus\{\zeta\}\subset\Omega$ and that $\zeta$ separates $E$ into two components $E^+$ and $E^-$ both tending to $\infty$. Consider $E^+$ to be the component tending to $\infty$ while nearing $L^+$ (see Figure \ref{parabolic-limsup}). Then, we denote by $H'$ the simply connected and convex in the positive direction domain bounded by $L^+, A$ and $E^-$.

\begin{figure}
    \centering
    \includegraphics[scale=0.5]{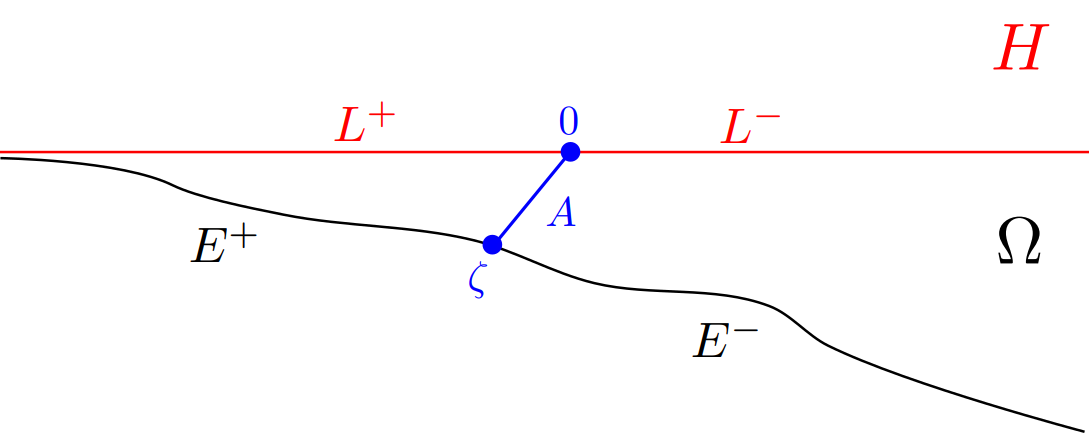}
    \caption{Sets of Theorem \ref{thm:asympttang}: \textbf{Case 3}}
    \label{parabolic-limsup}
\end{figure}

By a known formula (see \cite[p. 100]{ransford}), we have $\omega(x_{t_n},L^-,H)=\frac{\pi- \Arg x_{t_n}}{\pi}$ and thus
$$\lim\limits_{n\to+\infty}\omega(x_{t_n},L^-,H)=0.$$
Moreover, by the maximum principle for harmonic functions, we have $\omega(x_{t_n},L^-,H)\ge\omega(x_{t_n},E^-,H'),$ for all $n\in\mathbb{N}$. Hence $\omega(x_{t_n},E^-,H') \xrightarrow{n\to+\infty}0$.
By the Strong Markov Property, it follows that
$$\omega(x_{t_n},E^-,\Omega)=\omega(x_{t_n},E^-,H')+\int\limits_{L^+\cup A}\omega(w,E^-,\Omega)\cdot\omega(x_{t_n},dw,H'),$$
for all $n\in\mathbb{N}$. We need to evaluate the last integral and to this goal, we decompose it into two integrals. Firstly, 
$$\int\limits_{A}\omega(w,E^-,\Omega)\cdot\omega(x_{t_n},dw,H')\le\omega(x_{t_n},A,H').$$
However, the sequence $\{x_{t_n}\}$ converges to $\infty$, while $A$ is bounded. Directly, ${\omega(x_{t_n},A,H')\xrightarrow{n\to+\infty}0}.$
Modifying slightly the notation, we write
$$\int\limits_{L+}\omega(w,E^-,\Omega)\cdot\omega(x_{t_n},dw,H')=\int\limits_{(-\infty,0)}\omega(s,E^-,\Omega)\cdot\omega(x_{t_n},ds,H')=:I(n).$$ 
Since $E^+$ and $L^+$ tend to intersect, it is easily seen that $\omega(s,E^-,\Omega)\xrightarrow{s\to-\infty}0$. Let $\epsilon>0$. Then, there exists $s_0<0$ such that $\omega(s,E^-,\Omega)<\epsilon$ for all $s\le s_0$. As a result,
\begin{eqnarray*}
 I(n)&=&\int\limits_{(-\infty,s_0]}\omega(s,E^-,\Omega)\cdot\omega(x_{t_n},ds,H')+\int\limits_{(s_0,0)}\omega(s,E^-,\Omega)\cdot\omega(x_{t_n},ds,H')\\
 &\le&\epsilon+\omega(x_{t_n},(s_0,0),H').
\end{eqnarray*}
Once again, the fact that $(s_0,0)$ is a bounded set leads to $\omega(x_{t_n},(s_0,0),H')\xrightarrow{n\to+\infty}0$ and thus, $\limsup\limits_{n\to+\infty}I(n)\le\epsilon$. Nevertheless, this is true for all $\epsilon>0$, which means that $I(n)\xrightarrow{n\to+\infty}0.$
Combining everything and returning back to the Strong Markov Property, we deduce that 
$$\lim\limits_{n\to+\infty}\omega(x_{t_n},E^-,\Omega)=\lim\limits_{n\to+\infty}\omega(x_{t_n},E^-,H')=0.$$
In a similar fashion, we may prove that ${\omega(x_{t_n},F,\Omega)\xrightarrow{n\to+\infty}0}$, for any other component $F$ of the boundary $\partial\Omega$, other than $E^+$. Therefore, our construction necessitates that ${\omega(x_{t_n},E^+,\Omega)\xrightarrow{n\to+\infty}1}$. It can be easily understood that $E^+$ and $\partial\Omega\setminus E^+$ correspond through $h^{-1}$ to two circular arcs $B^+, B^-\subset\partial\mathbb{D}$, respectively, that intersect at two points, one of which is $\tau$. The conformal invariance of the harmonic measure certifies that 
$$\lim\limits_{n\to+\infty}\omega(h^{-1}(x_{t_n}),B^+,\mathbb{D})=1$$
and in view of Remark \ref{caratheodory}, we understand that $\{h^{-1}(x_{t_n})\}$ converges to $\tau$ tangentially. But the points $h^{-1}(x_{t_n}), n\in\mathbb{N}$, lie on the geodesic $\eta$ and therefore, the convergence should be orthogonal. Contradiction! Tracing back our initial hypothesis, we find that indeed $\liminf\limits_{t\to-\infty}|h(z)+t-x_t|>0$.

If alternatively $\partial H\cap\partial\Omega\ne\emptyset$, a similar proof, albeit considerably simpler, leads to the same result. In particular, in case $\partial H=\partial\Omega$, the result is straightforward.

Finally, a very similar procedure shows that the hypothesis that $\limsup\limits_{t\to-\infty}|h(z)+t-x_t|<+\infty$ leads to the same kind of contradiction. The main reason behind this is that in order for the limit superior to be finite, the real part of a sequence $\{x_{t_n}\}$ should converge to $-\infty$, while the imaginary part should remain bounded. Thus, we can resort to the same estimates as before to reach the aforementioned conclusion.

Summing up, we obtain $|h(z)+t-x_t| \xrightarrow{t\to-\infty}+\infty$ and
$$\limsup\limits_{t\to-\infty}\min\{\delta_\Omega(h(z)+t),\delta_\Omega(x_t)\}\le\lim\limits_{t\to-\infty}\delta_\Omega(h(z)+t)=c.$$
As a consequence, through Lemma \ref{lem:distancelemma}, $v_\Delta^T(t)\xrightarrow{t\to-\infty}+\infty$.
\end{proof}
\end{theorem}

\begin{remark}
    Combining the different cases of Theorem \ref{thm:asympttang} with the knowledge about backward orbits and petals, we understand that
     $\limsup\limits_{t\to-\infty}v_\Delta^T(t)<+\infty$ if and only if $\phi_t(z)$ converges non-tangentially, as $t\to-\infty$, to a point of the unit circle, for some (and equivalently all) $z\in\Delta$. Alternatively, $\limsup\limits_{t\to-\infty}v_\Delta^T(t)<+\infty$ if and only if $\phi_t(z)$ converges, as $t\to-\infty$, to a repelling fixed point of $(\phi_t)$, for some (and equivalently all) $z\in\Delta$.
\end{remark}

\begin{remark}
    The part of the proof of \textbf{Case 1} which correlates the ratio $\frac{\delta_\Omega^+(x_{t_n})}{\delta_\Omega^-(x_{t_n})}$ with the tangential convergence of the sequence is not an original work. It is a minor modification of the proof of \cite[Theorem 1.2]{quasi-backwards}. We include this part for the sake of clarity and completeness of the present article.
\end{remark}

Taking Corollary \ref{cor:distinparabpetal}, Theorem \ref{thm:divratehyperbolic} and Theorem \ref{thm:asympttang} along with Pythagoras' Theorem into account, the following result concerning the orthogonal and the tangential speed appears. We omit the proof since it is straightforward.
\begin{corollary}\label{cor:orthogonalasymb}
Let $(\phi_t)$ be a semigroup in $\mathbb{D}$ and $\Delta$ a petal of $(\phi_t)$. The following are true. 
\begin{enumerate}
    \item[\textup{(i)}] If $\Delta$ is hyperbolic with repelling spectral value $\lambda\in(-\infty,0)$, then 
    $$\lim\limits_{t\to-\infty}\frac{v^o_\Delta(t)}{t}=\frac{\lambda}{2} \quad \text{and} \quad \lim\limits_{t\to-\infty}\frac{v^T_\Delta(t)}{t}=0.$$
    \item[\textup{(ii)}] If $\Delta$ is parabolic, then
    $$\lim\limits_{t\to-\infty}\frac{v_\Delta^o(t)}{t}=\lim\limits_{t\to-\infty}\frac{v_\Delta^T(t)}{t}=0.$$
\end{enumerate}

\end{corollary}

\section{Non-Regular Backward Orbits}\label{sec:non-regular}

Throughout the current section, we analyze the asymptotic behavior of the speed of convergence along non-regular backward orbits. 
As discussed in Section \ref{semigroups}, a non-regular backward orbit for a one-parameter semigroup can either lie on the boundary of a petal, or converge to a super-repelling fixed point of the semigroup. 

Suppose $(\phi_t)$ is a semigroup of holomorphic self-maps of $\D$ and  $\gamma:[0,+\infty)\to\mathbb{D}$ is a non-regular backward orbit for $(\phi_t)$. 
Through the Koenigs function $h$, we move to the associated planar domain $\Omega$. Then the image $h(\gamma[0,+\infty))$ is either a half-line that converges to $\infty$ through the negative direction or a half-spiral that converges to $\infty$. The set $h(\gamma[0,+\infty))$ is either contained in a boundary component of a petal or has as a starting point a super-repelling fixed point. 

The following result indicates that the asymptotic behavior of the ``generalized'' total speed along non-regular backward orbits depends neither on the type of the semigroup nor the type of a petal, on whose boundary component the non-regular backward orbit may lie. 

\begin{proposition}\label{prop:non-regular}
There exists a one-parameter semigroup $(\varphi_t)_{t\geq0}$ of holomorphic self-maps of $\D$ such that $$ \lim_{t \to -\infty}\frac{d_{\D}(\zeta, \varphi_t(\zeta))}{t^2} =0, $$
where $\zeta$ lies on some non-regular backward orbit $\gamma$ of $(\varphi_t)$. 
Furthermore, there exists a one-parameter semigroup $(\psi_t)_{t\geq0}$ of holomorphic self-maps of $\D$ such that $$\liminf_{t \to -\infty} \frac{d_{\D}(\zeta, \psi_t(\zeta))}{t^2} \geq \frac{1}{4}, $$
where $\zeta$ lies on some non-regular backward orbit $\tilde{\gamma}$ of $(\psi_t)$. 
\begin{proof}
The proof is based on the existence of hyperbolic petals for non-elliptic semigroups of holomorphic self-maps of $\D$. Following the same argumentation, one can prove the same result working on the boundary of a parabolic petal or on a backward orbit converging to a super-repelling fixed point. Moreover, adjusting suitably the proof for spirallike domains and utilizing Lemma \ref{lem:mapofpetal}, the result implies generalization to the case of elliptic semigroups. 

Let $\Omega$ be a simply connected subdomain of $\mathbb{C}$ which is also convex in the positive direction and contains a maximal strip $S$. Let us denote by $\partial S^{+}$ and $\partial S^{-}$ the upper and lower boundary components of $S$, respectively. Suppose $\zeta \in \partial S^{+}$. We further assume that there exists some $t_0<0$ such that for all $t\leq t_0 $, ${\delta_{\Omega}(\zeta+t)= (\ln (-t))^{-1}}$. We note that due to maximality, it should be true that $\delta_{\Omega}(\zeta+t) \xrightarrow{t\to -\infty}0$. According to Lemma \ref{lem:distancelemma}, for some $t\leq t_0$ we obtain
\begin{eqnarray*}
\frac{d_{\Omega}(\zeta,\zeta+t)}{t^2} &\leq& \frac{d_{\Omega}(\zeta,\zeta+t_0)}{t^2} + \frac{d_{\Omega}(\zeta+t_0,\zeta+t)}{t^2} \\
&\leq &\frac{d_{\Omega}(\zeta,\zeta+t_0)}{t^2}  + \frac{1}{t^2} \int_{[\zeta+t,\zeta+t_0 ]}  \frac{ \di s}{\delta_{\Omega}(s )} \\
&=&\frac{d_{\Omega}(\zeta,\zeta+t_0)}{t^2} +  \frac{1}{t^2} \int_{-t_0}^{-t} \ln s \, \di s \\ 
&=&\frac{d_{\Omega}(\zeta,\zeta+t_0) +t_0 \ln(-t_0) -t_0}{t^2} + \frac{t-t \ln (-t) }{t^2} \xrightarrow{t\to -\infty}0.
\end{eqnarray*}
Let $h$ be the Riemann mapping of $\Omega$. We define the non-elliptic semigroup $(\varphi_t)_{t\geq 0}$ with ${\varphi_t(z) : = h^{-1}(h(z)+t)}$, for $z\in \D$, and $t\geq 0$. Then the maximal strip $S$ corresponds to a hyperbolic petal of $(\varphi_t)$ and $h^{-1}(\zeta) $ lies on a non-regular backward orbit. Thus, due to the conformal invariance of the hyperbolic distance, it follows that
$$\lim_{t \to -\infty}\frac{d_{\D}(h^{-1}(\zeta), \varphi_t(h^{-1}(\zeta)))}{t^2}=0.$$

For the second part, we again assume that $\Omega^{\prime}$ is a simply connected subdomain of $\mathbb{C}$ which is convex in the positive direction and contains a maximal strip $S^{\prime}$. 
Suppose again $\zeta \in \partial (S^{ \prime})^ {+}$. We further assume that there exists some $t_0<0$ such that for all $t\leq t_0 $, $\delta_{\Omega^{\prime}}(\zeta+t)= -  te^{-t^2}$. 
According to the Lemma \ref{lem:distancelemma}, for some $t\leq t_0$ we obtain
\begin{eqnarray*}
\frac{d_{\Omega^{\prime}}(\zeta,\zeta+t)}{t^2} &\geq&  \frac{d_{\Omega^{\prime}}(\zeta+t_0,\zeta+t)}{t^2} -\frac{d_{\Omega^{\prime}}(\zeta,\zeta+t_0)}{t^2}\\
&\geq &\frac{1}{4t^2} \log \left( 1+ e^{t^2} \frac{|t-t_0|}{-t}\right)- \frac{d_{\Omega^{\prime}}(\zeta,\zeta+t_0)}{t^2} \xrightarrow{t\to -\infty} \frac{1}{4}. 
\end{eqnarray*}
Applying the same technique as in the previous case, we can construct a non-elliptic semigroup $(\psi_t)_{t\geq 0}$, where $S^{\prime}$ is the image of a hyperbolic petal under the Koenigs function and $\zeta$ lies on the image of a non-regular backward orbit, with the use of the Riemann mapping $g$ of $\Omega^{\prime}$. Due to the conformal invariance of the hyperbolic distance, it follows 
$$\liminf_{t \to -\infty}\frac{d_{\D}(g^{-1}(\zeta), \psi_t(g^{-1}(\zeta)))}{t^2}\geq\frac{1}{4}.$$
\end{proof}
\end{proposition}

\begin{corollary}
    The asymptotic behavior of the speeds of convergence along non-regular backward orbits depends solely on the euclidean geometry of the Koenigs domain $\Omega$ and how fast the image of the non-regular backward orbit under the associated Koenigs function approaches asymptotically the boundary of $\Omega$.
\end{corollary}

\end{document}